\documentclass[11pt,leqno]{amsart}
\usepackage{latexsym,amssymb}
\usepackage{hyperref}
\usepackage{amsmath,amsthm}
\usepackage{amsfonts}
\usepackage[american]{babel}
\usepackage{mathrsfs}
\usepackage{psfrag}

\def\capa{{\rm Cap}_p}
\def\comp{\mathbf{R}^n\setminus}

\newcommand{\R}{\mathbb{R}}
\newcommand{\rn}{\mathbb{R}^N}

\numberwithin{equation}{section}
\newtheorem{theorem}{Theorem}[section]

\begin{document}
\title[Optimal Concavity of potentials]{Characterization of balls through optimal concavity for potential functions}

\author[P. Salani]{Paolo Salani}

\address{P. Salani, Dip.to di Matematica ``U. Dini'', Universit\`a degli Studi di Firenze, Viale Morgagni 67/A, 50134 Firenze - Italy}
\email{salani@math.unifi.it}
\date{}

\keywords{Capacity, power concave functions, Brunn-Minkowski inequality, convexity}
\subjclass{35R35, 35J65, 35B05}

\begin{abstract}
Let $p\in(1,n)$. If $\Omega$ is a convex domain in $\rn$ whose $p$-capacitary potential function $u$  is $(1-p)/(n-p)$-concave (i.e. $u^{(1-p)/(n-p)}$ is convex), 
then $\Omega$ is a ball.
\end{abstract}

\maketitle

\section{Introduction}
Let $n\geq 3$, $\Omega\subset\mathbf{R}^n$
and $p\in(1,n)$. The $p$--capacity of $\Omega$ can be defined
as follows (see for instance \cite{Evans-Gariepy}, \S 4.7):
\begin{equation}
\label{capa}
\capa(\Omega)=\inf\left\{\int_{\mathbf{R}^n}|\nabla v|^p\,dx\,:\,
v\in C^{\infty}_c(\mathbf{R}^n)\,,\, v(x)\ge 1\,\,\forall\,
x\in\Omega\right\}\,,
\end{equation}
where $C^{\infty}_c(\mathbf{R}^n)$ denotes the set of functions from
$C^{\infty}(\mathbf{R}^n)$ having compact support.

In the sequel $\Omega$ is a bounded open convex set, then
the above infimum is in fact a minimum which is realized by  the (classical) solution $u$
of the following problem
\begin{equation}
\label{pb1}
\left\{
\begin{array}{lll}
\textrm{div}(|\nabla u|^{p-2}\nabla u)=0\,,\quad\textrm{in}\quad\comp\overline\Omega\\
{}\\
u(x)=1\quad\textrm{in }\overline\Omega\,,\\
\\
\lim_{|x|\to+\infty}u(x)=0\,.
\end{array}
\right.
\end{equation}
The function $u$ is called \emph{$p$--capacitary potential function} of $\Omega$ and it holds
\begin{equation}\label{pcapaint}
\capa(\Omega)=\int_{\mathbf{R}^n\setminus\overline\Omega}|\nabla u|^p\,dx\,.
\end{equation}

It is well known that if $\Omega$ is (bounded, open and) convex, then $u$ is {\em quasi-concave}, that is all its superlevel sets
$$
\Omega(t)=\{x\in\R^n\,:\,u(x)\geq t\}\quad t\in(0,1]
$$
are convex, see \cite{gab, lew, kaw}. In fact, if $\Omega$ is smooth and strictly convex, one could even expect $u$ to satisfy some stronger concavity property, namely that there exist some suitable 
$\alpha(\Omega)<0$ such that $u$ is $\alpha(\Omega)$-concave, see Section \ref{alphaconcsec}. We recall here that a positive function is said {\em $\alpha$-concave}, for $\alpha<0$, if $u^{\alpha}$ is convex (see again Section \ref{alphaconcsec} for more details).

Indeed, when $\Omega$ is a ball of radius $R>0$ centered at $x_0$, it is easy to find explicitly the solution of \eqref{pb1}, that is
$$
u(x)=R^{q}\,|x-x_0|^{-q}
$$
where
$$
q=\frac{n-p}{p-1}\,,
$$
and it results to be $(-1/q)$-concave.

In this short note I prove that nothing better is possible and that this power concavity is optimal, in the sense that the property of {\em $u^{-1/q}$ to be convex} characterizes balls. Precisely, the main result of this paper is the following.
\begin{theorem}\label{mainthm}
Let $p\in(1,n)$, $\Omega$ be a bounded convex domain in $\rn$ and $u$ be the solution of \eqref{pb1}.
If $u$ is $(-1/q)$-concave, where $q=(n-p)/(p-1)$, then $\Omega$ is a ball.
\end{theorem}

To prove this theorem we will use three main ingredients:\\
- the first one is the Brunn-Minkowski inequality for $p$-capacity and its equality condition, proved in \cite{Borell, CJL} for $p=2$ and in \cite{CoSa1} for a generic $p$;\\
- the second ingredient is an easy relation existing between the $p$-capacity of a generic level set of $u$ and the capacity of $\Omega$, see formula \eqref{CapOmegat};\\
- the third ingredient is the expression of $p$-capacity through the behaviour at infinity of the potential function,
see formula \eqref{limuinfty}.
\medskip

In fact the last ingredient is needed to prove the following property, which may have its own interest and it is new, to my knowledge.

\begin{theorem}\label{secondthm}
If the solution $u$ of \eqref{pb1} has two homothetic level sets, then $\Omega$ is a ball.
\end{theorem}
In particular: {\em if $u$ has a level set that is homothetic to $\Omega$, then $\Omega$ is a ball}.
We recall here that two sets $A,B\subset\rn$ are said homothetic if there exist
$\rho>0$ and $\xi\in\rn$ such that $B=\rho A+\xi$, i.e. if they are dilate and translate of each other.
\medskip

To some extent, both the problems considered in Theorem \ref{mainthm} and Theorem \ref{secondthm} fall in the framework of overdetermined problems: in the first case the overdetermination is given by the concavity property of the solution $u$ of \eqref{pb1}, in the latter case the overdetermination is given by the existence of two homothetic level sets of $u$. 
\medskip

The paper is organized as follows. Firstly in Section 2 I introduce notation and recall some needed results and formulas (in particular the three main ingredients recalled above).  I prove Theorem \ref{secondthm} in Section 3.
Finally Section 4 contains the proof of Theorem \ref{mainthm}.

\section{Preliminaries}

\subsection{Basic notation.}
If $a,b\in\rn$, we denote by $\langle a,b\rangle$ their
scalar product and
by $|a|$ the euclidean norm of the vector $a$, i.e.
$|a|=\sqrt{\langle a,a\rangle}$.

If $M$ is an $n\times n$ symmetric matrix,
we denote by $\textrm{tr}(M)$ and
$\det(M)$ its trace and its determinant respectively;
$M > 0$  means that $M$ is positive definite, $M^{T}$ is the transposed of $M$ and $M^{-1}$ its						 inverse.

If $C$ is a subset of $\mathbb{R}^n$, $|C|$ is its Lebesgue measure, $\overline C$ is its closure, $\text{int}(C)$ is its interior and $\partial C$ is its boundary

For $r>0$ and $x\in\rn$, we denote by $B(x,r)$ the ball of radius $r$ centered at $x$. Then we set
$B^{N}=B(0,1)$ and $\omega_N=|B^{N}|$. By $S^{N-1}$ we denote the unit sphere in $\rn$, that is
$S^{N-1}=\partial B^{N}=\{x\in\rn\,:\,|x|=1\}$. Then $|S^{N-1}|=n\omega_N$.

If $u$ is a twice differentiable function, by $D u$ and $D^2 u$
we denote, as usual, the gradient of $u$ and its Hessian matrix
respectively, i.e.
$D u=(\frac{\partial u}{\partial x_1},\dots,
\frac{\partial u}{\partial x_N})$ and
$D^2 u=(\frac{\partial^2 u}{\partial x_i\partial x_j})_{i,j=1}^N$ and we denote by $||u||_{L^{p}(\Omega)}$ the $L^{p}$ norm of the function $u : \Omega \mapsto \mathbb{R}$.

\subsection{Ingredient 1: the Brunn-Minkowski inequality for $p$-capacity.}

The original form of the Brunn--Minkowski inequality involves volumes of
convex bodies (i.e. compact convex subsets of
$\mathbf{R}^n$ with non--empty interior) and states that $\textrm{Vol}_n(\cdot)^{1/n}$ is a concave function with respect to the
Minkowski addition, i.e.
\begin{equation}
\label{1II}
\left[\textrm{Vol}_n(\lambda K_1+(1-\lambda)K_2)\right]^{\frac{1}{n}}\ge
\lambda\,\left[\textrm{Vol}_n(K_1)\right]^{\frac{1}{n}}
+(1-\lambda)\,\left[\textrm{Vol}_n(K_2)\right]^{\frac{1}{n}}
\end{equation}
for every convex bodies $K_1$ and $K_2$ and $\lambda\in[0,1]$. Here
$\textrm{Vol}_n$ is the $n$-dimensional Lebesgue measure and the
Minkowski addition of convex sets is defined as follows
$$
A + B = \{ x + y \, \,|\, \, x \in A, \, \, y \in B\}\,,
$$
while $\lambda A=\{\lambda x\,:\,x\in A\}$ for any $\lambda\in\R$, as usual.

Inequality (\ref{1II}) is one of the fundamental results
in the modern theory of convex bodies; it can be extended to measurable
sets and several other important inequalities, e.g. the isoperimetric
inequality, can be deduced from it.

Suitable versions of the Brunn-Minkowski inequality hold also for the other quermassintegrals (see \cite{Schneider, Gardner}) and recently
Brunn-Minkowski type inequalities have been proved for several important geometric and analytic functionals
(see for instance \cite{Borell, bl, Colesanti, CoSa1, LiuMaXu, Salani1, Salani2} and especially the beautiful survey paper \cite{Gardner}).
Notice that in all the known cases, equality conditions are the same as in the classical Brunn-Minkowski inequality
for the volume, i.e. equality holds if and only if the involved sets
are (convex and) homothetic (i.e. translate and dilate of each other).

We will use the following theorem from \cite{CoSa1}.

\begin{theorem}[Theorem 1, \cite{CoSa1}]
\label{teorema1I} Let $K_1$ and $K_2$ be $n$--dimensional convex bodies and
let $p\in(1,n)$. Then
\begin{equation}
\label{1I}
\left[\capa(\lambda K_1+(1-\lambda)K_2)\right]^{\frac{1}{n-p}}\ge
\lambda\,\left[\capa(K_1)\right]^{\frac{1}{n-p}}+(1-\lambda)\,
\left[\capa(K_2)\right]^{\frac{1}{n-p}}\,,
\end{equation}
for every $\lambda\in[0,1]$. Moreover equality holds if and only if
$K_1$ and $K_2$ are homothetic.
\end{theorem}

Roughly speaking (\ref{1I}) says that $\capa(\cdot)^{\frac{1}{n-p}}$ is
a concave function in the class of convex bodies endowed with the Minkowsky
addition. But what is most relevant to the present paper is the equality condition: 
if equality holds in \eqref{1I}, then $K_1$ and $K_2$ are homothetic.

We recall that in the case of
the Newton capacity, i.e. for $p=2$ and $n\ge 3$, inequality (\ref{1I}) was proved by C. Borell \cite{Borell} and more recently in \cite{CJL} L.A. Caffarelli, D. Jerison and E.H. Lieb treated
the equality case. In \cite{CoSa1} the treatments of the inequality and of its equality case are unified and the results are extended to a generic $p\in(1,n)$.

\subsection{Ingredient 2: the $p$-capacity of a level set of the potential.}
Let $u$ be the $p$-capacitary potential of a domain $\Omega$ and set
$$
\Omega(t)=\{x\in\rn\,:\,u(x)\geq t\}\,
$$
for some $\tau\in(0,1]$.
Then it is easily seen that the following holds
\begin{equation}\label{CapOmegat}
\capa(\Omega(t))=t^{1-p}\capa(\Omega)\,.
\end{equation}
Indeed, the $p$-capacitary potential $u_t$ of $\Omega(t)$ is given by
$u_t(x)=t^{-1}u(x)$, as it can be trivially verified, and \eqref{CapOmegat} follows directly from
\eqref{pcapaint} or \eqref{limuinfty}.

\subsection{Ingredient 3: an expression of $p$-capacity through the behavior at infinity of the potential.}
In the case $p=2$ ($n\geq 3$) it is well known that the following relation between the Newton capacity of a convex domain and the behavior at infinity of the newtonian potential holds:
\begin{equation}\label{p2limuinfty}
{\rm Cap}_2(\Omega)=(n-2)\omega_n\lim_{|x|\to\infty}u(x)|x|^{n-2}\,.
\end{equation}
An analogous relation holds in the generic case $p\in(1,n)$

\begin{equation}\label{limuinfty}
\capa(\Omega)=c_{n,p}\left[\lim_{|x|\to\infty}u(x)|x|^{(n-p)/(p-1)}\right]^{p-1}\,,
\end{equation}
where
$$
c_{n,p}=\left(\frac{n-p}{p-1}\right)^{p-1}\omega_n\,;
$$
refer to \cite{CoSa1} for instance.

\subsection{Power-concave and quasi-concave functions.}\label{alphaconcsec}
For $\alpha\neq 0$ a function $v:\rn\to[0,+\infty]$ is said {\em $\alpha$-concave} if $\frac{\alpha v^\alpha}{|\alpha|}$ is concave.
In other words: if $\alpha>0$, $v$ is $\alpha$-concave if $v^\alpha$ is concave; if $\alpha<0$, $v$ is $\alpha$-concave if $v^\alpha$ is convex.
For $\alpha=0$, we say that $v$ is {\em log-concave} if $\log v$ is concave.

Furthermore, $v$ is said {\em quasi-concave} if all its super level sets $\{x\in\rn\,:\,v(x)\geq t\}$ are convex. To some extent, quasi-concavity corresponds to $\alpha$-concavity when $\alpha=-\infty$.

It is easily seen that if $v$ is $\alpha$-concave (for some $\alpha\in\R$), then $v$ is $\beta$-concave for every $\beta\leq\alpha$. 
Moreover, it is obvious that every $\alpha$-concave function (for some $\alpha\in\R$) is quasi-concave.

Given a quasi-concave function $v$, it is then natural to ask if it satisfies some better concavity properties and following \cite{ken} we define the 
{\em concavity number} of $u$ as follows
$$
\alpha(v)=\sup\{\beta\,:\,v\text{ is }\beta\text{-concave}\}\,.
$$
By \cite[Property 5]{ken}, for any $C^2$ quasi-concave function it is possible to explicitly calculate
$$
\alpha(v)=1-\sup\{v(x)v_{\theta\theta}(x)v_\theta(x)^{-2}\,:\,x\in\rn,\,\theta\in S^N,\,v_\theta(x)\neq 0\}\,,
$$
where $v_\theta$ and $v_{\theta\theta}$ denote respectively the first and the second derivatives of $v$ in direction $\theta$.

According to \cite{BLS, CoSa1, LS}  it is also possible and useful to associate to any quasi-concave function $v$ a {\em support function}
$h_v:\R\times\rn\to\R$, such that $h_v(X,t)$ is the support function of the super level set $\{v\geq t\}$ calculated at $X$, i.e.
$$
h_v(X,t)=\sup\{\langle X, x\rangle\,: v(x)\geq t\}\,.
$$ 
In this way the concavity of $v$ corresponds to the concavity of $h_v$ with respect to $t$, that is $v$ results to be concave if and only if
$$
\frac{\partial^2 h_v(X,t)}{\partial t^2}\leq 0\quad\text{for every }X\in\rn,\,t\in\R\,.
$$
Consequently $v$ is $\alpha$-concave (for some $\alpha\neq 0$) if and only if
$$
(1-\alpha)\frac{\partial h_v(X,t)}{\partial t}+t\,\frac{\partial^2 h_v(X,t)}{\partial t^2}\leq 0 \quad\text{for }(X,t)\in\rn\times\R\,.
$$
Notice that
$$
\frac{\partial h_v}{\partial t}\leq 0
$$
since $v$ is quasi-concave. Then we can also write
$$
\alpha(v)=\inf\left\{1+t\,\frac{\partial^2 h_v}{\partial t^2}\left(\frac{\partial h_v}{\partial t}\right)^{-1}\,:\,(X,t)\in\rn\times\R\right\}\,.
$$

As already said in the Introduction, it is well known that when $\Omega$ is convex, its $p$-capacitary potential $u$ is quasi-concave; then
we set
$$
\alpha(\Omega)=\alpha(u)\,.
$$
If $\Omega$ is sufficiently regular and strictly convex, one can expect that $\alpha(\Omega)>-\infty$. The aim of this paper is to prove the following:
{\em if $\alpha(\Omega)=(1-p)/(n-p)$, then $\Omega$ is necessarily a ball.}

\section{Proof of Theorem \ref{secondthm}}
Let $0<r<s\leq 1$, $\rho>1$, $\xi\in\rn$ such that
\begin{equation}\label{homolevel}
\Omega(r)=\rho\,\Omega(s)+\xi\,,
\end{equation}
that is $\Omega(r)$ and $\Omega(s)$ are the homothetic superlevel sets of the statement.

Notice that, since $r<s$, it holds 
$$\Omega(s)\subset\Omega(r)\,.
$$

For $t\in(0,1]$, let us denote by $u_t$ the $p$-capacitary potential of $\Omega(t)$, i.e. the solution of
$$
\left\{
\begin{array}{lll}
\textrm{div}(|\nabla u_t|^{p-2}\nabla u_t)=0\,,\quad\textrm{in}\quad\comp\Omega(t)\\
{}\\
u_t(x)=1\quad\textrm{in }\overline\Omega(t)\,,\\
\\
\lim_{|x|\to+\infty}u_t(x)=0\,.
\end{array}
\right.
$$
Then
$$
u_t(x)=\frac{u(x)}{t}\,,\quad x\in\comp\Omega(t)\,.
$$
In particular
$$
u_r(x)=\frac{u(x)}{r}\quad\text{for }x\in\comp\Omega(r)\quad\text{and}\quad
u_s(x)=\frac{u(x)}{s}\quad\text{for }x\in\comp\Omega(s)\,.
$$
On the other hand by \eqref{homolevel} it holds
$$
u_r(x)=u_s\big(\frac{x-\xi}{\rho}\big)
$$
and we finally get
\begin{equation}\label{urus}
u(x)=\frac{r}{s}\,u\big(\frac{x-\xi}{\rho}\big)\,,\quad x\in\comp\Omega(r)\,.
\end{equation}
Then by \eqref{limuinfty} the latter implies
$$
\begin{array}{rl}
\capa(\Omega)&=\lim_{|x|\to\infty}u(x)|x|^{\frac{n-p}{p-1}}=
\dfrac{r}{s}\lim_{|x|\to\infty}\,u\big(\dfrac{x-\xi}{\rho}\big)|x|^{\frac{n-p}{p-1}}\\
\\
&=\dfrac{r}{s}\,\rho^{\frac{n-p}{p-1}}\,\lim_{|x|\to\infty}\,u\big(\dfrac{x-\xi}{\rho}\big)
\left(\dfrac{|x-\xi|}{\rho}\right)^{\frac{n-p}{p-1}}\left(\dfrac{|x|}{|x-\xi|}\right)^{\frac{n-p}{p-1}}\\
\\
&=\dfrac{r}{s}\,\rho^{\frac{n-p}{p-1}}\,\capa(\Omega)\,,
\end{array}
$$
whence
\begin{equation}\label{rhors}
\frac{r}{s}=\rho^{\frac{p-n}{p-1}}\,.
\end{equation}
Moreover \eqref{urus} entails
\begin{equation}\label{nuovalabel}
\Omega(t)=\rho\,\Omega(\frac{s}{r}\,t)+\xi\quad\text{for }t\geq r\,.
\end{equation}
Hence, by setting
$$
s_k=\left(\frac{r}{s}\right)^ks=\rho^{\frac{k(p-n)}{p-1}}s\,,\quad k=0,1,\dots,
$$
it holds
$$
\lim_{k\to\infty}s_k=0
$$
and
$$
\Omega(s_k)=\rho\,\Omega(s_{k-1})+\xi=\rho^2\Omega(s_{k-2})+\rho\xi+\rho=\dots=
\rho^k\Omega(s_0)+\xi\,\sum_{i=0}^{k-1}\rho^i=\rho^k\Omega(s)+\xi\,\frac{\rho^k-1}{\rho-1}\,.
$$
Now let $x,y\in\partial\Omega(s)$, i.e.
$$
u(x)=u(y)=s\,,
$$
and set
$$
x_k=\rho^kx+\xi\,\frac{\rho^k-1}{\rho-1}\,,$$
$$
y_k=\rho^ky+\xi\,\frac{\rho^k-1}{\rho-1}\,.
$$
Then
$$
\lim_{k\to\infty}|x_k|=\lim_{k\to\infty}|y_k|=\infty
$$
and \eqref{limuinfty} yields
\begin{equation}\label{uxkuyk}
\lim_{k\to\infty}u(x_k)|x_k|^{\frac{p-1}{n-p}}=\capa(\Omega)=
\lim_{k\to\infty}u(y_k)|y_k|^{\frac{p-1}{n-p}}\,.
\end{equation}
On the other hand
$$
u(x_k)=u(y_k)=s_k\,,
$$
hence \eqref{uxkuyk} reads
$$
\lim_{k\to\infty}s_k|x_k|^{\frac{p-1}{n-p}}=
\lim_{k\to\infty}s_k|y_k|^{\frac{p-1}{n-p}}\,,
$$
that is
$$
\lim_{k\to\infty}\left(\rho^{\frac{k(p-n)}{p-1}}s\,\left|\rho^kx+\xi\,\frac{\rho^k-1}{\rho-1}\right|^{\frac{p-1}{n-p}}\right)=
\lim_{k\to\infty}\left(\rho^{\frac{k(p-n)}{p-1}}s\,\left|\rho^ky+\xi\,\frac{\rho^k-1}{\rho-1}\right|^{\frac{p-1}{n-p}}\right)\,,
$$
which implies
$$
\lim_{k\to\infty}\,\left|x+\xi\,\frac{1-\rho^{-k}}{\rho-1}\right|^{\frac{p-1}{n-p}}=
\lim_{k\to\infty}\,\left|y+\xi\,\frac{1-\rho^{-k}}{\rho-1}\right|^{\frac{p-1}{n-p}}\,.
$$
Since $\rho>1$, the latter finally implies
$$
\left|x+\xi\,\frac{1}{\rho-1}\right|=\left|y+\xi\,\frac{1}{\rho-1}\right|=R\,,
$$
which means that $\Omega(s)$ is a ball or radius $R$ centered at the point $\xi/(1-\rho)$,
i.e.
$$
\Omega(s)=B(\frac{\xi}{1-\rho},R)\,,
$$
and from \eqref{nuovalabel} we easily obtain
$$
\Omega(s_k)=B(\frac{\xi}{1-\rho},R\rho^k )\quad\text{for }k=1,2,\dots\,.
$$
Then $u$ is radial in $\rn\setminus\overline{\Omega(s)}$ and, by analytic continuation, it is radial
in $\rn\setminus\overline\Omega$ and $\Omega$ is a ball.

\section{Proof of Theorem \ref{mainthm}}

Let $u\in C(\rn)\cap C^2(\rn\setminus\overline\Omega)$ be the solution of \eqref{pb1} and 
$$q=(1-p)/(n-p)\,.
$$ 

We will proceed by proving that, if $u$ is $q$-concave, then all its level sets are homothetic.
Then the proof will be concluded thanks to Theorem \ref{secondthm}.
\medskip

Assume that
$$
v=u^q\quad\text{ is convex in }\rn\,.
$$
Hence for every $v_0,v_1\in\R$ and for every $\lambda\in(0,1)$ it holds
\begin{equation}\label{convexityofv}
\{x\in\rn\,:\, v(x)\leq (1-\lambda)v_0+\lambda v_1\}\supseteq(1-\lambda)\{x\in\rn\,:v(x)\leq v_0\}+
\lambda\,\{x\in\rn\,:\,v(x)\leq v_1\}\,.
\end{equation}
Now take $r,s\in(0,1]$, fix $\lambda\in(0,1)$ and set
\begin{equation}\label{t}
t=[(1-\lambda)r^q+\lambda\,s^q]^{1/q}\,.
\end{equation}
By setting $v_0=r^q$ and $v_1=s^q$, we have
$$
t^q=(1-\lambda)v_0+\lambda\,v_1
$$
and 
$$
\Omega(r)=\{v\leq r^q\}\,,
$$
$$
\Omega(s)=\{v\leq s^q\}\,,
$$
$$
\Omega(t)=\{v\leq t^q\}\,.
$$
Then \eqref{convexityofv} entails
$$
\Omega(t)\supseteq (1-\lambda)\,\Omega(r)+\lambda\,\Omega(s)\,.
$$
Thanks to the monotonicity of $p$-capacity with respect to set inclusion, the latter implies
$$\capa(\Omega(t))\geq\capa\big((1-\lambda)\,\Omega(r)+\lambda\,\Omega(s)\big)$$
and by chaining with the Brunn-Minkowski inequality for $p$-capacity (see Theorem \ref{teorema1I}) we obtain
\begin{equation}\label{usingbm}
\capa(\Omega(t))^{1/(n-p)}\geq (1-\lambda)\,\capa(\Omega(r))^{1/(n-p)}+\lambda\,\capa(\Omega(s))^{1/(n-p)}
\end{equation}

On the other hand, by \eqref{CapOmegat} we have
$$\capa(\Omega(r))=r^{1/(n-p)}\capa(\Omega)\,,
$$
$$
\capa(\Omega(s))=s^{1/(n-p)}\capa(\Omega)
$$
and
$$
\capa(\Omega(t))=t^{1/(n-p)}\capa(\Omega)\,.
$$
Substituting in \eqref{usingbm}  and taking into account \eqref{t},
we finally get 
$$
\capa(\Omega(t))^{1/(n-p)}=(1-\lambda)\,\capa(\Omega(r))^{1/(n-p)}+\lambda\,\capa(\Omega(s))^{1/(n-p)}\,,
$$
i.e. equality holds in \eqref{usingbm}, and consequently equality must hold in the Brunn-Minkowski inequality for $p$-capacity for $\Omega(r)$ and $\Omega(s)$. Then Theorem \ref{teorema1I} tells that
$\Omega(r)$ and $\Omega(s)$ must be homothetic. This concludes the proof thanks to Theorem \ref{secondthm}, as already said.

%
%
%
%
%

\end{document}